\documentclass{ifacconf}

\usepackage{xcolor}
 \usepackage{amsmath, amssymb}
\usepackage{graphicx}      
\usepackage{float}

\usepackage{filecontents}

\makeatletter
\let\old@ssect\@ssect 
\makeatother

\usepackage{natbib}
\usepackage{hyperref}

\makeatletter
\def\@ssect#1#2#3#4#5#6{%
  \NR@gettitle{#6}
  \old@ssect{#1}{#2}{#3}{#4}{#5}{#6}
}
\makeatother

\begin{document}
\begin{frontmatter}

\title{Towards Model Discovery Using Domain Decomposition and PINNs }


\author[TUBS]{Tirtho S. Saha}
\author[TUD]{Alexander Heinlein} 
\author[TUBS]{Cordula Reisch}

\address[TUBS]{Institute for Partial Differential Equations, TU Braunschweig, Germany (e-mail: c.reisch@tu-bs.de, t.saha@tu-bs.de)}
\address[TUD]{Delft Institute of Applied Mathematics, Delft University of Technology, The Netherlands (e-mail: a.heinlein@tudelft.nl).}

\begin{abstract}               
We enhance machine learning algorithms for learning model parameters in complex systems represented by ordinary differential equations (ODEs) with domain decomposition methods. The study evaluates the performance of two approaches, namely (vanilla) Physics-Informed Neural Networks (PINNs) and Finite Basis Physics-Informed Neural Networks (FBPINNs), in learning the dynamics of test models with a quasi-stationary longtime behavior. We test the  approaches for data sets in different dynamical regions and with varying noise level. 
As results, we find a better performance for the FBPINN approach compared to the vanilla PINN approach, even in cases with data from only a quasi-stationary time domain with few dynamics.  
\end{abstract}


\end{frontmatter}

\section{Introduction}
 Mathematical modeling of biological processes is inherently complex due to the intricate and often only partially understood mechanisms involved. Additionally, biological processes exhibit different behaviors on different temporal and spatial scales. 
 Some processes may take a long time, and often data are available only from certain stages, while data for other stages are unavailable. Using numerical methods to solve these inverse problems can be cumbersome. To alleviate this issue, alternative methods have emerged as a prominent solution. Furthermore, advanced approaches like domain decomposition-based Physics-Informed Neural Networks (PINNs) are gaining prominence in solving the model problem. This study addresses these challenges by leveraging data-driven machine learning approaches to identify abstract mechanisms and determine parameter values that represent the known biological mechanisms from observed data. 
 
 We focus on two toy models: the saturated growth model, which captures population growth dynamics, and the competition model, which examines interactions between two species, including scenarios of coexistence and survival (\cite{murray2007mathematical}). These models were tested on synthetic data from different time intervals, such as a dynamical and stationary phase, and the total time domain, with varying noise levels. 
 We employ two different approaches:
 physics-informed neural networks (PINNs, \cite{raissi2017physics}) using the SciANN library (\cite{haghighat2021sciann}), 
 and domain decomposition-based PINNs using finite basis PINNs (FBPINNs, \cite{moseley2023finite}); for the application of FBPINNs to (systems) of ordinary differential equations, we also refer to~\cite{heinlein_multifidelity_2024}.
 The aim is to compare the ability of the methods in learning the parameters of the dynamical system in cases where the data is limited to certain time intervals. 
 Problems like this occur when dealing with (biological) problems where only stationary data is available that can be interpreted as the result of a dynamical process in advance of the measurement.

 Up to our knowledge, the application of the domain decomposition approach is new in the field of parameter estimation, in particular for differential equations with data from quasi-stationary dynamics.

 \section{Computational Methods}\label{sec:computational}
 
 We start with introducing the two computational methods, first vanilla PINNs and then the more sophisticated  FBPINNs including the idea of domain decomposition.

 \subsection{Physics Informed Neural Networks(PINNs)}
 In contrast to purely data-driven approaches, PINNs are trained by using a combination of labeled training data and available prior knowledge about the problem (\cite{raissi2017physics},\cite{lagaris1998artificial},\cite{dissanayake_neural-network-based_1994}). In the forward problem setup, the physical law(s) are known and encoded in a PDE, but the solution of the PDE is unknown. Let us consider a differential equation of the general form:
\begin{equation}\label{eq:GeneralPDE}
    \begin{aligned}
        \mathcal{D}[u(x)] &= f(x),  \quad x \in \Omega \subset \mathbb{R}^d, \\
        \mathcal{B}_k[u(x)] &= g_k(x), \quad x \in \Gamma_k \subset \partial \Omega,
    \end{aligned}
\end{equation}
where $\mathcal{D}[u(x)] $ is some differential operator with $u(x)$ as the solution, and $\mathcal{B}_k[\cdot]$ is a boundary operator including as well the initial conditions (ICs), which ensure uniqueness of the solution. The input $x$ could be spatial and/or temporal, where $d$ is the dimension of the domain. Equation~(\ref{eq:GeneralPDE}) can represent many differential equation problems, including linear and nonlinear problems, ODEs and PDEs, time-dependent and time-independent problems, and problems with an initial value, Dirichlet and Neumann boundary problems. 

To solve the differential equation~(\ref{eq:GeneralPDE}), PINNs use an NN to directly approximate the solution, i.e., 
\begin{equation}\label{eq:NNAnsatz}
\begin{aligned}
    u^\mathrm{PINN}(x; \mathbf{\theta}) \approx u(x),
\end{aligned}
\end{equation}
where $x$ is the input to the network and $\theta$ are the trainable parameters of the NN model.
The proposed general loss function from \cite{raissi2017physics} to train the PINNs model combines two influences,  
\begin{equation}\label{eq:GeneralLossFunc}
\begin{aligned}
    \mathcal{L}(\mathbf{\theta}) = \mathcal{L}_\mathrm{PDE}(\mathbf{\theta}) + \mathcal{L}_\mathrm{BC}(\mathbf{\theta}).
\end{aligned}
\end{equation}
The PDE based loss function is 
\begin{equation}\label{eq:PDELossFunc}
\begin{aligned}
    \mathcal{L}_\mathrm{PDE}(\mathbf{\theta}) &= \frac{\lambda_{\mathrm{phy}}}{N_I} \sum_{i=1}^{N_I} (\mathcal{D}[u^\mathrm{PINN}(x_i; \theta)] - f(x_i))^2  
\end{aligned}
\end{equation}
with  $\{x_i\}_{\substack{i=1}}^{N_I}$ is a set of collocation points sampled within $\Omega$ and $\lambda_{\mathrm{phy}}$ is a weight. 
The boundary condition loss function is
\begin{equation}\label{eq:BCLossFunc}
\begin{aligned}
    \mathcal{L}_\mathrm{BC}(\mathbf{\theta}) &= \sum_{k=1}^{N_k} \frac{\lambda_{B}^k}{N_{B}^k} \sum_{i=1}^{N_{B}^k} (\mathcal{B}_k[u^\mathrm{PINN}(x_{i}^{k}; \theta)] - g_k(x_{i}^{k}) )^2,
\end{aligned}
\end{equation}
where $\{x_{j}^k\}_{\substack{j=1}}^{N_{B}^{k}}$ is a set of points sampled along each boundary condition (BC) and $\lambda_{B}^k$ is a weight.
The weights $\lambda_I$ and $\lambda_{B}^k$ are chosen so that the individual terms in the loss function \eqref{eq:GeneralLossFunc} contribute in a balanced manner.
Finding an appropriate choice of $\lambda_I$ and $\lambda_{B}^k$ leading to the best result is usually challenging and problem-depending. 

An alternative to using separate boundary condition loss terms $\mathcal{L}_\mathrm{BC}(\mathbf{\theta})$ is to hard constrain the solution to satisfy the boundary condition exactly, which we do in this work. 
This approach involves directly incorporating the boundary conditions into the neural network architecture to inherently satisfy the boundary condition, thereby removing the need for the boundary condition residual term in the loss function.  
The hard constraining modifies the NN ansatz in~\eqref{eq:NNAnsatz} as following: 
\begin{equation}\label{eq:NNAnsatzConstrained}
\begin{aligned}
    \tilde{u}^\mathrm{PINN}(x; \mathbf{\theta}) = \mathcal{C}[u^\mathrm{PINN}(x; \mathbf{\theta})] \approx u(x),
\end{aligned}
\end{equation}
where $\mathcal{C}$ is the constraining operator applied to the output of the NN model. Consequently, the loss function of the NN becomes: 
\begin{equation}\label{eq:GeneralLossFuncConstrained}
\begin{aligned}
    \mathcal{L}(\mathbf{\theta}) = 
    \mathcal{L}_\mathrm{PDE}(\mathbf{\theta}) = \frac{\lambda_{\mathrm{phy}}}{N_I} \sum_{i=1}^{N_I} (\mathcal{D}[\tilde{u}^\mathrm{PINN}(x_i; \theta)] - f(x_i))^2  .
\end{aligned}
\end{equation}

In many real-world scenarios, the objective is not only to solve a forward problem but also to address an inverse problem. An inverse problem involves estimating unknown parameters or initial conditions based on observed data, with the governing equations explicitly defined or partially known. In this case, the differential and/or boundary operators $\mathcal{D}$ respectively $\mathcal{B}_k$ may depend on a set of additional parameters $P = (p_1,\ldots,p_{N_P})$. Hence, solving the inverse problem does not only involve finding the network parameters $\theta$ but also the parameters $P$. Given available synthetic or real data and other prior knowledge, the inverse problem reads
\begin{equation}\label{eq:GeneralLossFuncWithData}
\begin{aligned}
    \min_{\theta,P} \mathcal{L}(\mathbf{\theta},P) &= \mathcal{L}_\mathrm{PDE}(\mathbf{\theta},P)  + \mathcal{L}_\mathrm{data}(\mathbf{\theta},P) +  \mathcal{L}_\mathrm{par}(\mathbf{\theta},P),
\end{aligned}
\end{equation}
where the loss function in~\eqref{eq:GeneralLossFuncConstrained} of the PINNs approach is complemented by the data loss
\begin{equation}\label{eq:dataLossFunc}
\begin{aligned}
    \mathcal{L}_\mathrm{data}(\mathbf{\theta},P) &= \frac{\lambda_\mathrm{data}}{N_D} \sum_{i=1}^{N_D} (\tilde{u}^\mathrm{PINN}(x_{i}; \theta) - u_{\text{data}}(x_{i}))^2 
\end{aligned}
\end{equation}
and an optional parameter-induced loss function 
 \begin{equation}\label{eq:paraLossFunc}
\begin{aligned}
    \mathcal{L}_\mathrm{par}(\mathbf{\theta},P) &= \lambda_\mathrm{param} \sum_{i=1}^{N_P} (\max \{ 0, p_{i, \min} - p_i, p_i - p_{i,\max} \})^2.
\end{aligned}
\end{equation}
Here, $\lambda_\mathrm{data}$ represents the weight for the data residual in the loss function, and $N_D$ is the number of data points used for the training. In general, the data points $u_{\mathrm{data}}(x_i)$ could consist of measurements, observations, or synthetic data derived from simulations. Additionally, $\lambda_\mathrm{param}$ is the weight for the parameter constraints term in the loss function, and $N_P$ is the number of parameters. $(p_{i, \min}, p_{i, \max})$ are the \texttt{min} and \texttt{max} the values of the $i$'th parameter. The range for the values of $p_{i,\min}$ and $p_{i,\max}$ depends on known values, and in the specific cases of this work, all parameters are considered non-negative, so the lower threshold is known. 

As more terms are included in the loss function, the complexity of the training increases due to the higher number of interactions between terms and the additional constraints imposed on the optimization process.

\subsection{Domain Decomposition-Based PINNs(FBPINNs)}
Vanilla PINN approaches show a spectral bias, meaning that they can learn the low frequency components more easily than the high-frequency components of the solution. 
To address this issue, it has been observed that domain decomposition-based PINN architectures can learn multiscale components of the solution by \cite{moseley2023finite}. 

The domain decomposition-based PINN approach defines an approximate solution similar to those given in ~\eqref{eq:NNAnsatz} and ~\eqref{eq:NNAnsatzConstrained}, as it works well with both soft and hard boundary constraints setups. However, it differs in terms of the network architecture, as it employs as many neural networks as the number of subdomains chosen. The global network produces the output
\begin{equation}\label{eq:FBPINNsNetworkOutput}
    u^\mathrm{FBPINN}(x; \theta) = \sum_{j=1}^{J} \omega_j(x) \cdot \text{unnorm} \circ u^\mathrm{sub}_j (x; \theta_j) \circ \text{norm}_j(x),
\end{equation}
where the term $u^\mathrm{FBPINN}(x; \theta)$ represents the collective sum of the output of all subdomains. The normalization term $\text{norm}_j$ adjusts the input variable $x$ to the range of $[-1,1]$ in each dimension over the subdomain before it is input to the individual neural network $u^\mathrm{sub}_j(x; \theta_j)$ at $\Omega_j$. Then comes the output unnormalisation $\text{unnorm}$ term, which ensures that the output stays within the range $[-1, 1]$ in each neural network. Finally, the outputs are multiplied by the window function $\omega_j(x)$, which is smooth, differentiable, and zero outside the subdomain, confines the network's solution locally.
Moreover, the choice of the subdomains enhances the learning of specific frequencies fitting to the subdomain size. 

\begin{figure}[htb]
  \centering
  \includegraphics[width=7.4cm]{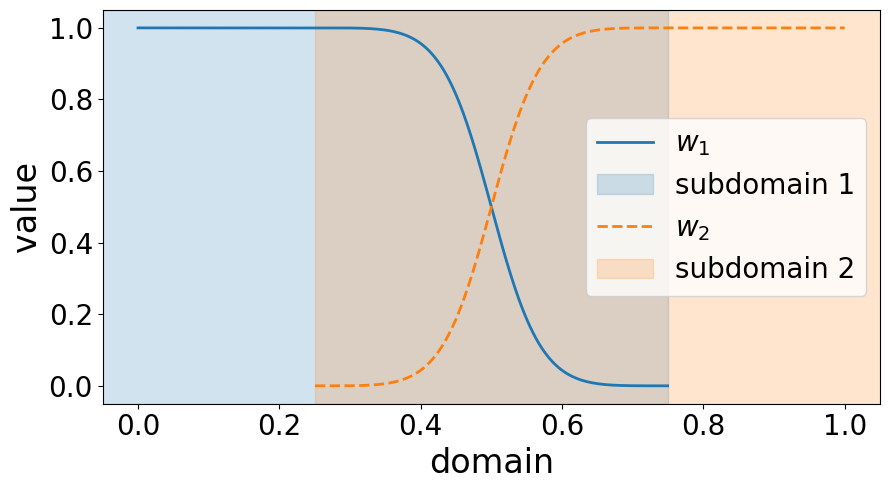} 
  \caption{Window function~\eqref{eq:SubdomainWindowFunction}. The light blue and light orange regions represent the respective subdomain intervals, while the combined light brown region highlights the overlap between the two subdomains.}
  \label{fig:WindowFunction}
\end{figure}

The window function depicted in Fig.~\ref{fig:WindowFunction}  for a hyperrectangular subdomain is defined from a partition of unity as
\begin{equation}\label{eq:SumOverWFto1}
        \sum_{j=1}^{J} \omega_j \equiv 1 \quad \text{on } \Omega,  \quad \quad \mathrm{supp}(\omega_j) \subset \overline{\Omega}_j,
\end{equation}
with
\begin{equation}\label{eq:SubdomainWindowFunction}
    \begin{aligned}
    \omega_j(x) &=  \left(H(x - x_{j,\min}) \cdot  H(x_{j,\max} - x)\right) \cdot \\& \frac{1}{4} \left(1 + \cos\left(\pi\frac{x - \mu_{j}}{\sigma_{j}}\right)\right)^2,
    \end{aligned}
\end{equation}
where $\omega_j$ is the window function in the $j$'th subdomain after domain decomposition, and $H$ is the Heaviside step function that ensures the solution is zero outside of each subdomain. 
The interval $(x_{j,\min}, x_{j,\max})$ represents the left and right overlapping region for the subdomain $j$. $\mu_{j}$ and $\sigma_{j}$ represent the center and half-width of each subdomain, respectively. 
The cosine function and the Heaviside step function ensures the solution is smooth within the interval, having a continuous first derivative, and zero outside of the subdomain overlapping region.

The loss function for the FBPINN method is calculated for the given ansatz defined in~\eqref{eq:NNAnsatzConstrained}, with the network output provided in~\eqref{eq:FBPINNsNetworkOutput}, similar to the PINN loss function in~(\ref{eq:GeneralLossFuncWithData}) giving
\begin{equation}\label{eq:FBPINNLossFuncWithDataAndPDE}
\begin{aligned}
    \mathcal{L}(\mathbf{\theta}) &= \mathcal{L}_\mathrm{PDE}(\mathbf{\theta})+\mathcal{L}_\mathrm{data}(\mathbf{\theta}) +\mathcal{L}_\mathrm{par}(\mathbf{\theta}),
\end{aligned}
\end{equation}
with the loss functions in \eqref{eq:GeneralLossFuncConstrained}, \eqref{eq:dataLossFunc} and \eqref{eq:paraLossFunc}, where the PINN solution $\tilde{u}^\mathrm{PINN}$ is replaced by the FBPINN solution~\eqref{eq:FBPINNsNetworkOutput}.

These two computational approaches will be used in the following for determining in an inverse problem setting the parameters of ordinary differential equations with available data in certain time domains. The domain decomposition approach therefore applies to the time domain, giving more weight to time domains with more data or, in perspective, time domains with higher quality data.

 \section{Mathematical Models}
We introduce two differential equation models with a non-trivial solution behavior. 
Those models will be investigated with the introduced methods of Sec.~\ref{sec:computational}. 
 
 \subsection{Saturated growth model}
\begin{figure}[H]
  \centering
  \includegraphics[width=8.4cm]{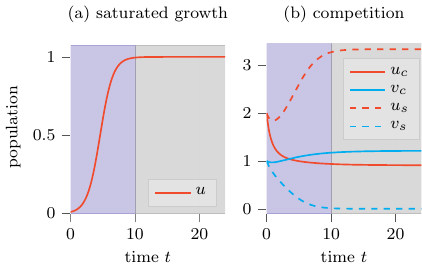}
  \caption{Solutions of (a) the saturated growth model \eqref{eq:SaturatedGrowthModel}  and (b) the competition model~\eqref{eq:CompetitionModel} for parameters with coexistence ($u_c, v_c)$ or single-survival $(u_s, v_s)$. The dynamic time frame is shaded in blue, the quasi-stationary in gray.}
  \label{fig:ODEModels}
\end{figure}

In the saturated growth model (e.g. \cite{murray2007mathematical}), we consider a population of one species $u$ with the carrying capacity $C$. The saturated growth model is
\begin{equation}\label{eq:SaturatedGrowthModel}
    \frac{du}{dt} = u(C - u) \quad \text{with} \quad u(0)=u_0,
\end{equation}
where $ u$ represents the population size, and $ C>0$ is the carrying capacity. This model captures the dynamics of a population undergoing saturated growth, such as the growth of a virus population in liver tissue. The solution tends towards $C$ for initial values $ u_0 > 0 $, representing the saturation of growth as the population reaches its maximum capacity. The growth rate is moderated by the term $(C-u)$, which implies that the population growth rate decreases as it approaches the carrying capacity $C$. Fig.~\ref{fig:ODEModels}(a) shows the solution of the saturated growth model for an initial value $ u(0) = u_0 > 0 $ and $C=1$.

\subsection{Competition Model}
The second model is a Lotka-Volterra competition model (e.g. \cite{murray2007mathematical}), which served as a test model in \cite{reisch2023model} concerning the possibility of model discovery. 
The dynamics of two species with inter- and inner-species competition is given by
\begin{equation}\label{eq:CompetitionModel}
    \begin{aligned}
        \frac{du}{dt} &= u(1 - a_1 u - a_2 v) =f_1(u,v), \\
        \frac{dv}{dt} &= rv(1 - b_1 u - b_2 v) =f_2(u,v), \\
        \text{with} \quad u(0)&=u_0>0; \quad v(0)=v_0>0,
    \end{aligned}
\end{equation}
where $a_1$, $a_2$, $b_1$, $b_2$ and $r$ are all positive coefficients, and the species $u$ and $v$ compete for shared resources. 
There are two non-trivial long-time behaviors depending on the parameter values. In the first case, the system tends toward a coexistence steady state $(u_c^*,v_c^*) \neq (0,0)$. In the second case, the system favors one species to steady state with other species vanishes, i.e., either $(u_s^*\neq 0,v_s^*=0)$ or $(u_s^*=0,v_s^* \neq 0)$. 

Fig. \ref{fig:ODEModels}(b) illustrates the competition models featuring coexistence and single-survival scenarios. As depicted, in coexistence scenarios, two species survive together within shared environments. Conversely, in survival scenarios, one species gradually dominates the other over time, determined by its initial value within the closed shared environments.
The parameter values used in the simulations are given in Tab.~\ref{tab:parameters_competition}.

\begin{table}
    \centering
     \caption{Parameter values of the competition model~\eqref{eq:CompetitionModel} in a coexistence and a single-survival setting. }
     \setlength{\tabcolsep}{12pt}
    \begin{tabular}{cccccc}
        parameter & $r$  &$a_1$  & $a_2$  &$b_1$  &$b_2$ \\ \hline
       coexistence  & 0.5 & 0.7 & 0.3 & 0.3 & 0.6\\
       single-survival  & 0.5 & 0.3 & 0.6 & 0.7 & 0.3\\
    \end{tabular}
    \label{tab:parameters_competition}
\end{table}

\subsection{Model discovery and parameter estimation}

The three models, saturated growth for one species and competition between two species with two parameter settings, serve as test models for a model discovery or, by restricting the mechanisms beforehand, for a parameter estimation in varying data scenarios.
The introduced PINN approaches solve the inverse problem of finding parameter values by including the parameter to learn in the underlying differential equation that contributes in $\mathcal{L}_\mathrm{PDE}$. In this study, we provide the terms in the ordinary differential equation that are included in the model used for data generation. In future works, we plan to implement a sparse choice of some mechanism terms, like in \cite{brunton2016discovering}. 

The parameters we want to learn here are $C$ for the saturated growth model and for the competition model the parameters in Tab.~\ref{tab:parameters_competition}.
Fig.~\ref{fig:ODEModels} shows in gray the different time intervals used for data sampling. 
We want to compare the ability of determining the parameters of the  models by using data either from the dynamical time interval $[0,10]$, the quasi-stationary time interval $[10,24]$, or in the whole time interval $[0,24]$. 
We will investigate both approaches, vanilla PINN and FBPINN, in all data settings. 

Our expectations for the computational results are based on analytical properties of the stationary points: 
In the saturated growth model~\eqref{eq:SaturatedGrowthModel}, the nontrivial stationary point $u^\star =C$ gives directly the parameter that we want to estimate. 
Our first hypothesis therefore is that both computational approaches should be able to reproduce a good estimation of the parameter, independent of the time frame employed for data collection. Besides, the solution of the neural network should be more precise when giving data from the dynamical time domain rather than from the quasi-stationary because the dynamics lead to the stationary states, but there are multiple dynamics tending towards the same stationary state.

The identification of the parameters in the competition model is much more challenging, firstly because there are more parameters to identify, and secondly because the stationary states do not depend on all parameters. 
More precisely, the coexistence stationary state is given by $$(u_c^\star, v_c^\star)=\left(\frac{a_2-b_2}{a_2b_1-a_1b_2}, \frac{a_1-b_1}{a_2b_1-a_1b_2}\right),$$ so it is independent of $r$ and we have in numerical simulations two stationary state values for determining four dependent parameters. 
 The stationary state in the single-survival parameter setting is $(u_s^\star, v_s^\star)=(1/a_1, 0)$, and hence, independent of $a_2, b_1, b_2, r$. By knowing only the longtime behavior of the solution, it is therefore hard to determine the parameter values except for $a_1$ in this setting. 

Our hypotheses for the competition model take this into account: 
We expect that learning $a_1$ in the single-survival setting is feasible, even when only datum in the quasi-stationary domain is available. On the other hand, determining any parameter in the coexistence case is hard for taking only data in the quasi-stationary domain. 

So far, these hypotheses on learning the parameters in the different models and settings are valid for both methods of Sec.~\ref{sec:computational}. 
We investigate now how the domain decomposition with overlapping domains affect (i) the parameter estimation and (ii) the quality of the learned NN solution, both compared to vanilla PINNs.

\section{Results}

We start with testing our hypotheses on the general parameter estimation problem and then compare more detailed the outcomes of vanilla PINNs and FBPINNs. 

\subsection{Model accuracy and parameter learning}

Our first hypothesis is that the parameter $C$ in the saturated growth model is easy to learn for any of the algorithms and independent of the data region used. 
This hypothesis is confirmed by test cases that we run with the hyperparameters in Tab.~\ref{tb:training_parameters}. 
The learned parameters $C$ for the three temporal domains and the two approaches, vanilla PINN and FBPINN, are given in Tab.~\ref{tab:saturated growth}.

\begin{table}
    \centering
        \caption{Learned values $C$ for the saturated growth model~\eqref{eq:SaturatedGrowthModel}, the true value is $C=1$.}
    \setlength{\tabcolsep}{18pt}    
    \begin{tabular}{cccc}
         & $[0,24]$ & $[0,10]$ & $[10,24]$ \\ \hline
      PINN    & 1.0042 & 0.9916 & 1.0097\\
      FBPINN   & 0.9917 & 0.9915 & 0.9995\\
    \end{tabular}

    \label{tab:saturated growth}
\end{table}

\begin{figure}[htb]
  \centering
  \includegraphics[width=8.4cm]{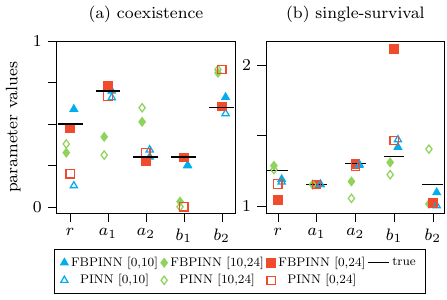}
  \caption{Learned parameters in the competition model. The values $[a,b]$ give the time domain of data used.}
  \label{fig:SciANN_FBPINN_LearnedParams}
\end{figure}

Next, we check our hypotheses on the competition model. 
One states that in the coexistence case, determining the parameters from the quasi-stationary time domain is impossible due to an underdetermined algebraic system for the parameters depending on the steady states. 
Fig.~\ref{fig:SciANN_FBPINN_LearnedParams} shows these problems with false estimates in the quasi-stationary setting but acceptable estimates in the dynamical or full time domain. Both approaches perform qualitatively same in this case.

The second hypothesis for the competition model states that in the single-survival parameter setting the estimation of $a_1$ is possible even in the quasi-stationary time domain, while it is not possible to determine the other parameters in this time domain exactly. 
Fig.~\ref{fig:SciANN_FBPINN_LearnedParams} supports this hypothesis, and again, both methods perform qualitatively equal. 
Consequently, the parameter estimation is a task that does not improve by FBPINNs, which was expected from the nature of the problem. 

\begin{figure}[H]
  \centering
  \includegraphics[width=8.4cm]{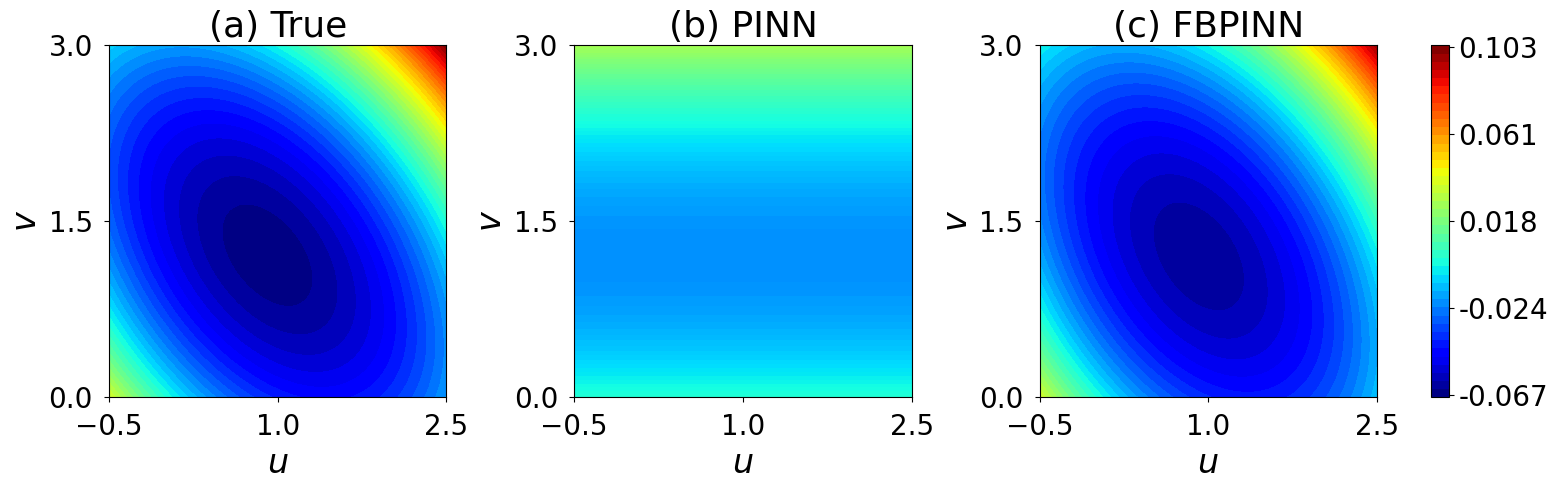}
  \caption{Energy plots for the competition model~\eqref{eq:CompetitionModel} in the coexistence setting with data in $[0, 24]$.}
  \label{fig:EnergyPlots}
\end{figure}

 However, the effect of learned parameters from different approaches can be observed in Fig.~\ref{fig:EnergyPlots} by the energy plots using the Lyapunov function 
  \begin{equation}\label{eq:LyapunovFuncCompModel}
    \begin{aligned}
        \phi =& -a_1 b_2 r (b_1 u + a_2 v) + a_1 a_2 b_1 b_2 r u v \\&+ \frac{1}{2} r a_1 b_2 (a_1 b_1 u^2 + a_2 b_2 v^2).
    \end{aligned}
\end{equation}

 Even though the differences of the learned parameters in the time interval $[0,24]$ are relatively small, the energy plots given by the Lyapunov function based on the learned parameters differ rather crucial. 
 While the learned parameters with FBPINN give an energy functional that strongly resembles the energy functional of the ground truth parameters, the vanilla PINN energy functional is qualitatively different. 
 This has an effect on variations of the initial conditions: While the solutions will still tend towards a point close to the true stationary state in the FBPINN case, the dynamics with the learned parameters from vanilla PINN may be totally different.

\begin{figure}[htb]
  \centering
  \includegraphics[width=8.4cm]{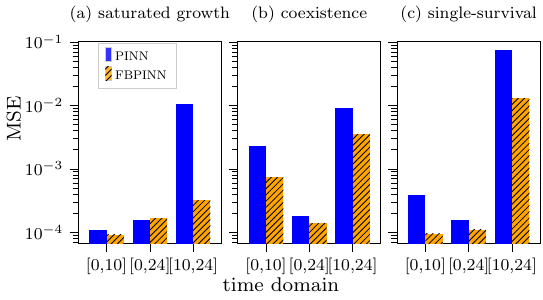}
    \caption{MSE of the vanilla PINN solution and the  FBPINN solution based on data from three time intervals. }
  \label{fig:MSEVanillaPINNsVsFBPINNs}
\end{figure}

Next, we want to compare the second outcome of the PINN approaches, a learned solution of the NN. 
Fig.~\ref{fig:MSEVanillaPINNsVsFBPINNs} shows the mean squared error of the NN output and the numerical solution of~\eqref{eq:SaturatedGrowthModel}, resp.~\eqref{eq:CompetitionModel}, as ground truth. 
The results of the comparison in Fig.~\ref{fig:MSEVanillaPINNsVsFBPINNs} show a higher accuracy of the FBPINN compared to vanilla PINN. This difference is very prominent as well for the quasi-stationary time domains, where the parameter estimation failed in both parameter settings of the competition model.

\begin{figure}[H]
  \centering
  \includegraphics[width=8.4cm]{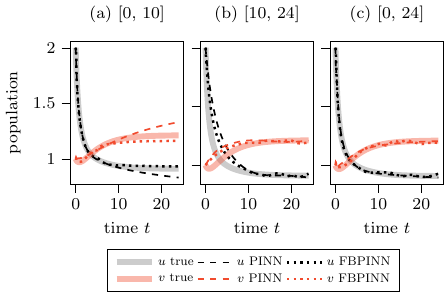}
  \caption{Comparison of time-dependent PINN and FBPINN solutions  for the competition model with coexistence}
    \label{fig:CoexistenceModelComparison}
\end{figure}

Based on this impression, we can dive deeper into the differences of the solutions for one case.
The difference of the MSE is largest for the competition models. Therefore, we compare the time-resolved solutions of the competition model in the coexistence case, see Fig.~\ref{fig:CoexistenceModelComparison}. 
The differences in the models depend, of course, on the data time domains. 
The solution of the vanilla PINN in the dynamical time domain $[0,10]$ has a larger error for larger time, while the error for the solution with quasi-stationary data from $[10,24]$ has a larger error for small time. 
The FBPINN solution in the whole time domain shows a surprising oscillatory behavior for larger time. 
A reason for this behavior might be overfitting of the included noise. Therefore, the oscillations only occur when data in the quasi-stationary time domain is available.
A more sophisticated hyperparameter tuning may reduce this effect.

To our surprise, the overlap size of the subdomains had only a small influence on the solution quality with one exception: In the single survival test case, we found a large MSE for a window overlap into the data region of $wo=1.001$ and smaller. 
For the other test cases, a window overlap variation between $1.001$ and $2.3$ does not affect the solution quality. 
In all regarded cases, there were still collocation points in the overlap. 

\subsection{Loss landscapes with and without domain decomposition}

\begin{figure}[htb]
  \centering
  \includegraphics[width=8.4cm]{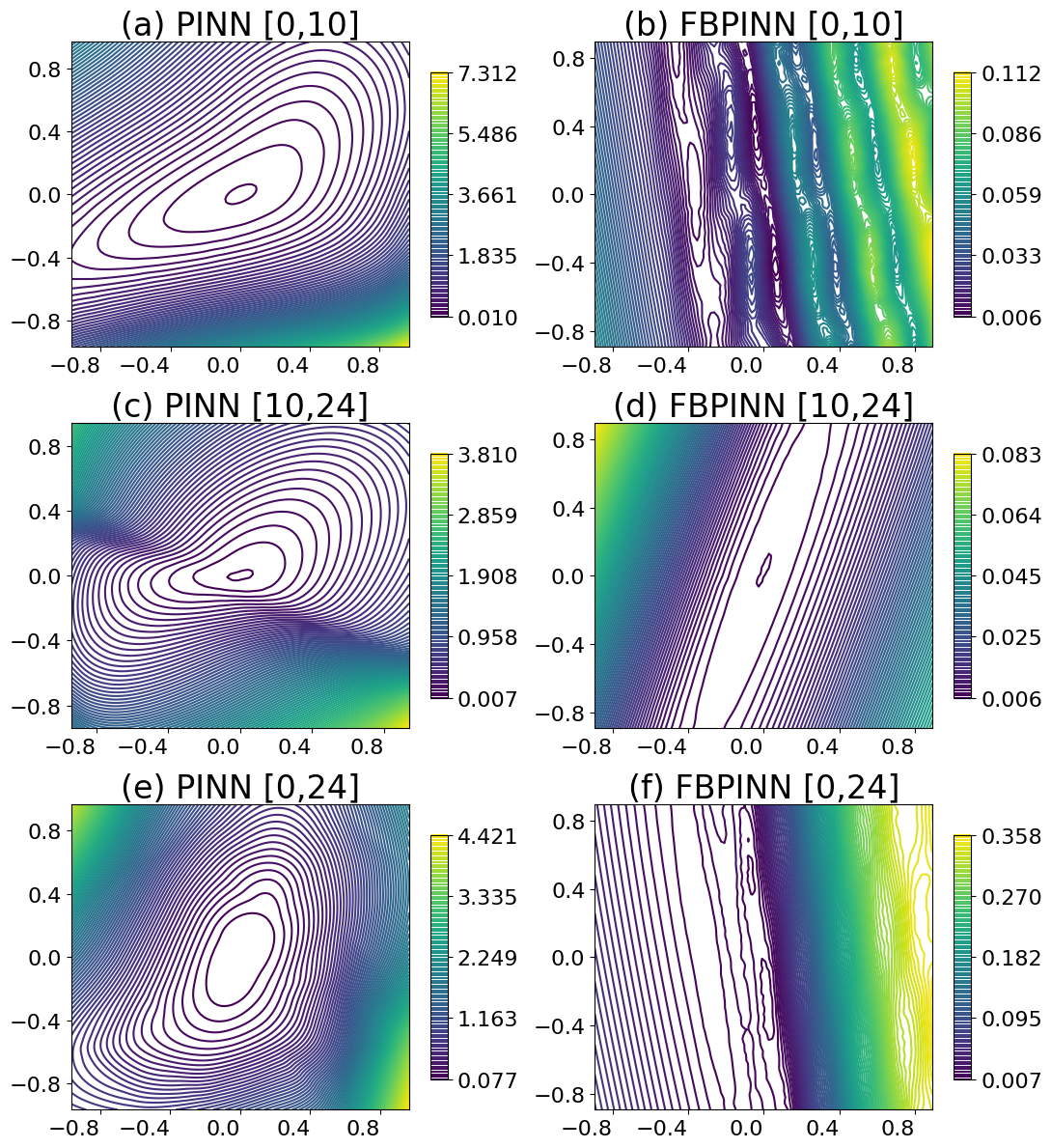}
    \caption{PINN  and FBPINN loss landscapes for the competition model with coexistence in the three data settings.}
  \label{fig:LossLandscape}
\end{figure}

The loss landscapes (\cite{li2018visualizing}) of FBPINNs and PINNs are shown in Fig.~\ref{fig:LossLandscape} with individual colorbars.
In total, the loss landscape of FBPINNs shows less sensitivity to small variations in the trained weights compared to those of vanilla PINNs.
The PINN loss landscapes are more convex than the loss landscapes for FBPINNs in the data regions $[0,10]$ and $[0,24]$. 
Following \cite{li2018visualizing}, this may indicate that the network initialization is more crucial for FBPINNs, where some chaotic regions exist next to well-formed convex regions. 
Further interpretation of the loss landscapes is challenging because the landscape shows only two random directions of the large parameter space of the NNs.

\subsection{Noise Effect}

\begin{figure}[htb]
      \centering
      \includegraphics[width=7.4cm]{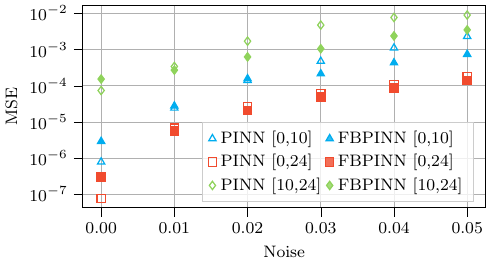} 
            \caption{MSE of solutions for the competition model in the coexistence case with varying noise for both, the PINN and the FBPINN approach.}
      \label{fig:MSECoexistenceModel}
    \end{figure}

The time-dependent dynamics of the learned solutions in Fig.~\ref{fig:CoexistenceModelComparison} shows some oscillatory behavior due to noise. Fig.~\ref{fig:MSECoexistenceModel} compares the MSE for different noise levels in the data for the competition model with coexistence. 

    For any non-zero noise level, the FBPINN solution has a smaller MSE than the vanilla PINN solution. In some cases, the difference is in the order of one magnitude. 
    This supports the better performance of FBPINNs compared to PINNs without any domain decomposition. 
    A reason for this observation may be the ability of FBPINN to learn different parts of the solution more stable due to the different subdomains. 

\section{Training parameters}

The used hyperparameters are given in Tab.~\ref{tb:training_parameters}.
The results are stable against changes of $nC$, $wo$, $wi$ and the number of layers.
The FBPINNs approach uses all the training points in a single training step. Consequently, we set the batch size to be equal to the sum of the number of data points $(nD)$ and the number of collocation points $(nC)$, for the PINNs approach as well. 
The code is available at \\ \href{https://github.com/tirtho109/VanillaPINNsVsFBPINNs}{github.com/tirtho109/VanillaPINNsVsFBPINNs}.

\begin{table}[H]
\begin{center}
\caption{Hyperparameters }
\label{tb:training_parameters}
\begin{tabular}{ccc}
\textbf{Parameter} & \textbf{PINN} & \textbf{FBPINN} \\\hline
Hidden Layers (layers) & [5,5,5] & [5,5,5] \\
Epochs & 50000 & 50000 \\
Activation Function & tanh & tanh \\
Physics loss weight $(\lambda_{\text{phy}})$ & 1.0 & 1.0 \\
Data loss weights $(\lambda_{\text{data}})$ & 1.0 & 1.0 \\
Optimizer & Adam & Adam \\
Learning rate & 0.001 & 0.001 \\
Number of collocation points ($nC$) & 200 & 200 \\
Collocation Sampling & Grid & Grid \\
Number of data points ($nD$) & 100 & 100 \\
Batch Size & 300 & 300 \\
Number of subdomains ($nsub$) & N/A & 2 \\
Window overlap data-region ($wo$) & N/A & 1.9 \\
Window overlap no-data-region ($wi$) & N/A & 1.0005 \\
Parameter loss weight $(\lambda_{\text{param}})$ & N/A & $1 \times 10^6$ \\\hline
\end{tabular}
\end{center}
\end{table}

\bibliography{ifacconf}         

\begin{thebibliography}{5}
\providecommand{\natexlab}[1]{#1}
\providecommand{\url}[1]{\texttt{#1}}
\providecommand{\urlprefix}{URL }
\expandafter\ifx\csname urlstyle\endcsname\relax
  \providecommand{\doi}[1]{doi:\discretionary{}{}{}#1}\else
  \providecommand{\doi}{doi:\discretionary{}{}{}\begingroup
  \urlstyle{rm}\Url}\fi

\bibitem[{Brunton et~al.(2016)Brunton, Proctor, and
  Kutz}]{brunton2016discovering}
Brunton, S.L., Proctor, J.L., and Kutz, J.N. (2016).
\newblock Discovering governing equations from data by sparse identification of
  nonlinear dynamical systems.
\newblock \emph{Proceedings of the national academy of sciences}, 113(15),
  3932--3937.

\bibitem[{Dissanayake and
  Phan-Thien(1994)}]{dissanayake_neural-network-based_1994}
Dissanayake, M.W.M.G. and Phan-Thien, N. (1994).
\newblock Neural-network-based approximations for solving partial differential
  equations.
\newblock \emph{Communications in Numerical Methods in Engineering}, 10(3),
  195--201.

\bibitem[{Haghighat and Juanes(2021)}]{haghighat2021sciann}
Haghighat, E. and Juanes, R. (2021).
\newblock {SciANN}: A keras/tensorflow wrapper for scientific computations and
  physics-informed deep learning using artificial neural networks.
\newblock \emph{Computer Methods in Applied Mechanics and Engineering}, 373,
  113552.

\bibitem[{Heinlein et~al.(2024)Heinlein, Howard, Beecroft, and
  Stinis}]{heinlein_multifidelity_2024}
Heinlein, A., Howard, A.A., Beecroft, D., and Stinis, P. (2024).
\newblock Multifidelity domain decomposition-based physics-informed neural
  networks and operators for time-dependent problems.
\newblock ArXiv:2401.07888 [cs, math].

\bibitem[{Lagaris et~al.(1998)Lagaris, Likas, and
  Fotiadis}]{lagaris1998artificial}
Lagaris, I.E., Likas, A., and Fotiadis, D.I. (1998).
\newblock Artificial neural networks for solving ordinary and partial
  differential equations.
\newblock \emph{IEEE transactions on neural networks}, 9(5), 987--1000.

\bibitem[{Li et~al.(2018)Li, Xu, Taylor, Studer, and
  Goldstein}]{li2018visualizing}
Li, H., Xu, Z., Taylor, G., Studer, C., and Goldstein, T. (2018).
\newblock Visualizing the loss landscape of neural nets.
\newblock \emph{Advances in neural information processing systems}, 31.

\bibitem[{Moseley et~al.(2023)Moseley, Markham, and
  Nissen-Meyer}]{moseley2023finite}
Moseley, B., Markham, A., and Nissen-Meyer, T. (2023).
\newblock Finite basis physics-informed neural networks ({FBPINNs}): a scalable
  domain decomposition approach for solving differential equations.
\newblock \emph{Advances in Computational Mathematics}, 49(4), 62.

\bibitem[{Murray(2007)}]{murray2007mathematical}
Murray, J.D. (2007).
\newblock \emph{Mathematical biology: I. An introduction}, volume~17.
\newblock Springer Science \& Business Media.

\bibitem[{Raissi et~al.(2019)Raissi, Perdikaris, and
  Karniadakis}]{raissi2017physics}
Raissi, M., Perdikaris, P., and Karniadakis, G.E. (2019).
\newblock Physics-informed neural networks: A deep learning framework for
  solving forward and inverse problems involving nonlinear partial differential
  equations.
\newblock \emph{Journal of Computational Physics}, 378, 686--707.

\bibitem[{Reisch and Burmester(2023)}]{reisch2023model}
Reisch, C. and Burmester, H. (2023).
\newblock Model selection focusing on longtime behavior of differential
  equations.
\newblock ArXiv:2312.05128 [math].

\end{thebibliography}

\end{document}